\documentclass[12pt]{article}    %[birkhoff.tex   20.08.17]
\usepackage{amsfonts,amssymb} %,mlb}

\usepackage[margin=25mm, a4paper]{geometry}  % journal quality

\usepackage{tikz} 

\def\n{\noindent}  \def\?#1{} 
\def\IZ{{\mathbb{Z}}} \def\IN{{\mathbb{N}}} \def\IR{{\mathbb{R}}}
   \def\supp{{\rm supp}}
\def\ep{\varepsilon}    \def\cB{{\cal B}} \def\cM{{\cal M}} 
\def\cT{{\cal T}}  \def\t{\tilde} \def\bdelete#1{}
\def\beq#1#2{\begin{equation} \label{#1} #2 \end{equation}}
\def\bea#1{\begin{eqnarray*} #1 \end{eqnarray*}} \def\a{\!\!\!&\!\!\!\!&}
\def\toas#1{\stackrel{#1}{\longrightarrow}}
\def\function#1{\left\{\!\!\!\begin{array}{ll} #1 \end{array} \right.}
\def\proof{\smallskip \noindent {\bf Proof. \ }}       %start of proof
\def\blanksquare{\,\,\,$\sqcup\!\!\!\!\sqcap$}         %blank  square
\def\qed{\hfill\blanksquare\linebreak\smallskip\par}   %end of proof

\def\thname{Theorem}  \def\lmname{Lemma}    \def\prname{Proposition}
\def\dfname{Definition}  \def\crname{Corollary}  \def\rmname{Remark}
\def\exname{Example}

\newtheorem{theorem}{\thname}[section]   %Numbering: Theorem--Other section
\newtheorem{lemma}{\lmname}[section]     %{lemma}[theorem]{Lemma}   subsection
\newtheorem{proposition}[lemma]{\prname} %lemma
\newtheorem{corollary}[lemma]{\crname}   %lemma
\newtheorem{example}{\exname}[section]  

\newtheorem{dftn}{\dfname}[section]
\newenvironment{definition}{\begin{dftn}\rm}{\end{dftn}} %section
\newtheorem{rmrk}[lemma]{\rmname}
\newenvironment{remark}{\begin{rmrk}\rm}{\end{rmrk}}     %lemma

             \catcode`@=11 \@addtoreset{equation}{section} \catcode`@=12 

\makeatletter\def\fps@figure{htbp}\makeatother %figure pos: tbp - standard

\begin{document}

\title{Ergodic averaging with and without \\ invariant measures}
\author{Michael Blank\thanks{
        Institute for Information Transmission Problems RAS 
        (Kharkevich Institute);} 
        \thanks{National Research University ``Higher School of Economics''; 
        e-mail: blank@iitp.ru}
        \thanks{This work has been carried out at IITP RAS and we gratefully 
                acknowledge the support of Russian Foundation for Sciences 
                (project No. 14-50-00150).}
       }
\date{August 20, 2017} % \today
\maketitle

\begin{abstract}%
The classical Birkhoff ergodic theorem in its most popular version says 
that the time average along a single typical trajectory of a dynamical 
system is equal to the space average with respect to the ergodic 
invariant distribution. This result is one of the cornerstones of the 
entire ergodic theory and its numerous applications. Two questions 
related to this subject will be addressed: how large is the set of 
typical trajectories, in particular in the case when there are no 
invariant distributions, and how the answer is connected to properties 
of the so called natural measures (limits of images of ``good'' 
measures under the action of the system).
\end{abstract}%

%%%%%%%%%%%%%%%%%%%%%%%%%%%%%%%%
\section{Introduction}\label{s:intro}

A conventional {\em point-wise ergodic} {``Theorem''} says that: 
For a given invariant measure $\mu$ an ergodic average along a 
trajectory starting from a certain initial point converges $\mu$-almost 
everywhere.  
Examples: celebrated Birkhoff, Kingman, Oseledets theorems. 
Those results are well known but have two serious disadvantages. 
First, the ergodic invariant measure may have a small support in 
which case the behavior of trajectories starting from points outside 
of the support is not described by such statements. Second, as we 
will see there are simple examples of low dimensional dynamical 
systems having no invariant measures and thus formally having 
nothing to do with such claims. Our aim is to try to overcome these 
difficulties, namely to find conditions under which claims of ergodic 
averaging type may be obtained for reasonably general dynamical 
systems for almost all initial points with respect to a ``good'' reference 
measure (e.g. Lebesgue measure). 

Let $(X,\rho)$ be a compact metric space equipped with a $\sigma$-algebra of 
measurable sets $\cB$ and a probability reference measure $m$, and let 
$T:X\to X$ be a measurable map from this space into itself. In this paper we 
restrict ourselves to the questions related to the generalization of one of the 
most known and widely used results in ergodic theory of dynamical systems -- 
the classical Birkhoff ergodic theorem, which claims that % 
\beq{e:trj-conv}{\frac1n\sum_{k=0}^{n-1}f(T^kx) \toas{n\to\infty} \int f~d\mu }% 
for an ergodic $T$-invariant measure $\mu$, each integrable function $f\in L_\mu^1$ 
and $\mu$-a.e. $x\in X$. 

In a number of cases there is a special invariant measure called Sinai-Ruelle-Bowen (SRB) 
measure $\t\mu$ (exact definitions will be given in Section~\ref{s:basic}) which represents 
averages along trajectories starting from a set of positive $m$-measure. 
For a given probability measure $\mu$ consider the set of $\mu$-{\em typical points} %
\beq{e:star}{ Z_{\mu}:=\{x\in X:~~\frac1n\sum_{k=0}^{n-1}f(T^kx) \toas{n\to\infty} 
                                              \int f~d\mu \} \qquad \forall f\in C^0(X) \}, }% 
where $C^0(X)$ stands for the set of continuous functions on $X$. Our aim is to 
find conditions under which $Z_{\t\mu}$ is the set of full $m$-measure, i.e. $m(Z_{\t\mu})=1$. 
Moreover, we will show that the set of typical points may be large not only for ergodic 
invariant measures (which is not surprising), but also for some non-invariant measures 
(see Section~\ref{s:irr}). 

To start with, let us demonstrate that the situation when $m(\supp(\t\mu))$ is much smaller 
than $m(Z_{\t\mu})$ appears very naturally in the simplest examples of dynamical systems. 

\begin{example} \label{ex:1}
Let $T$ be a map from the unit disc $X:=\{(\phi,R):~0\le\phi<2\pi,~0\le R\le1\}$
into itself defined in the polar coordinates $(\phi,R)$ by the relation: %
\beq{e:circle}{T(\phi,R):=(\phi+2\pi\alpha+\beta(R-r)~{\rm mod~2\pi},~\gamma(R-r)+r)}  %
with parameters $\alpha,\beta,\gamma,r\in(0,1)$. 
\end{example}
The circle $\{R=r\}$ is the only attractor of this map with the basin of attraction 
consisting of all the points in $X$ except the unstable fixed point located at the origin. 
See Section~\ref{s:example} for the detailed discussion of ergodic properties of this system.

This example in fact is based on a trivial observation that a uniformly contractive 
dynamical system governed by the map $Tx:=x/2$ on the unit interval is uniquely 
ergodic (the Dirac measure $\delta_0$ at the origin is invariant) while 
$m(\supp(\delta_0))=0\ll 1=m(Z_{\delta_0})$.

Without a kind of attraction property for points outside of the support 
of the measure $\t\mu$ one cannot expect the generalization of the Birkhoff 
theorem, e.g. a presence of another attractor of positive reference measure 
obviously would contradict it. As we will see the contraction alone is also 
not enough for this.

In what follows we give necessary and sufficient conditions for the property 
$m(Z_{\t\mu})=1$ in three different cases: regular 
(Theorem~\ref{t:nsc} in Section~\ref{s:inv+}), 
irregular  (Theorem~\ref{t:suf-irr} in Section~\ref{s:irr}) and self-consistent 
(Theorem~\ref{t:self-con} in Section~\ref{s:self}). 
Definitions of these cases will be given in the corresponding sections. 
Roughly speaking the regular case corresponds to the situation when the map $T$ 
is smooth enough, in the irregular case due to discontinuities of the map 
$T$ there are no invariant measures, while the self-consistent case is a 
deterministic version of a nonlinear Markov chain. 
Technically the hard part of the proofs of all subsequent results is described is 
Section~\ref{s:proofs} which is devoted to the regular case only, while the 
proofs of corresponding claims made for more complicated situations discussed 
in Sections~\ref{s:irr} and \ref{s:self} 
are reduced to the those in Section~\ref{s:proofs} by means of a special trick: 
introduction of the notion of weakly ergodic measures. In Section~\ref{s:basic} 
we discuss basic definitions and constructions as well as some known results 
about ergodic averaging.

\section{Basic definitions and constructions}\label{s:basic}
We start with basic definitions and constructions used throughout the paper. 

\begin{definition} Let $(X,\cB)$ be a compact measurable space and let 
$T:X\to X$ be a measurable map from this space into itself. The map induces 
the {\em transfer-operator} $T_*$ acting in the space of probability measures 
$\cM$ on $X$ according to the formula $T_*\mu(A):=\mu(T^{-1}A)$ for 
each measurable set $A\subseteq X$. A measure $\mu\in\cM$ is called 
$T$-invariant if $T_*\mu=\mu$.  
\end{definition}

Throughout the paper we always assume that $\cB$ is a Borel $\sigma$-algebra 
of measurable sets and only probability measures will be taken into consideration. 
Additionally we fix a certain measure $m\in\cM$ on this $\sigma$-algebra being 
positive on each open measurable subset and call it a {\em reference} measure. 

\begin{definition} A measurable map $T:X\to X$ is called {\em nonsingular} with 
respect to the reference measure $m$ if for any $A\in\cB$, $T_*m(A)=0$ 
if and only if $m(A)=0$. $T$ is {\em conservative} if for each set $A\in\cB$ of 
positive $m$-measure there exists $n\in\IN$ such that $m(A\cap T^{-n}A)>0$. 
A measurable set is called {\em wandering}, if all its images under the action 
of the map $T$ are disjoint. 
A probability measure $\mu$ is called {\em ergodic} if $\mu(A)\in\{0,1\}$ 
for each $T$-invariant set. Recall that a set $A\subset\cB$ is $T$-{\em invariant} 
if $T^{-1}A=A$.   
\end{definition} 

Observe that the ergodicity of a probability measure does not necessarily imply 
its invariance, for example for each point $x\in X$ the $\delta$-measure $\delta_x$ 
is ergodic for any map $T:X\to X$.

\begin{definition} For a given probability measure $\nu$ denote by 
$\cM(\nu)\subseteq\cM$ the set of probability measures $\mu$ on $X$ absolutely 
continuous with respect to $\nu$ (notation $\mu\ll\nu$). 
\end{definition}

Recall that $\mu\ll \nu$ means that $\nu(A)=0$ implies $\mu(A)=0$. 
The convergence of measures throughout the paper is considered always in the weak* sense. 

\begin{definition} A measure $\mu$ is said to be {\em wandering} if all its images 
under the action of the transfer-operator $T_*$ are mutually singular. 
\end{definition}

Recall that the mutual singularity of measures $\mu, \nu$  (notation $\mu\perp\nu$) 
means that there exists a set $A\in\cB$ such that $\mu(A)=\nu(X\setminus A)=0$. 
Observe also that the presence of an wandering measure $\mu\in\cM(\nu)$ is 
equivalent the existence of an wandering set of positive $\nu$ measure.

\begin{definition} Fix a reference measure $m\in\cM$. A measure $\t\mu^{\rm nat}\in\cM$ 
is called {\em natural} if there is an open set $U\subseteq X$ with $m(U)>0$ 
such that %
\beq{e:nat}{ \frac1n\sum_{k=0}^{n-1}T_*^k\mu\toas{n\to\infty}\t\mu^{\rm nat}  
                   \qquad\forall \mu\in\cM(m), \quad \mu(U)=1 .}
A measure $\t\mu^{\rm obs}\in\cM$ is called {\em observable} if there is an open set 
$U\subseteq X$ with $m(U)>0$ such that %
\beq{e:obs}{\frac1n\sum_{k=0}^{n-1}T_*^k\delta_x\toas{n\to\infty}\t\mu^{\rm obs} 
                   \qquad\forall x\in U',} %
where $\delta_x$ is the Dirac measure at point $x$ and $m(U\setminus U')=0$. 
The set $U$ in these constructions is the {\em basin of attraction} of the corresponding 
measure. 
\end{definition}

The natural and observable measures are different instances of the so called  
Sinai-Ruelle-Bowen (SRB) measures. For the discussion of the connections 
between them we refer a reader to \cite{BB}, where the question when these 
objects coincide has been raised in the first time, and to \cite{Mi,JT} where further 
clarifications were obtained. A number of nontrivial examples demonstrating not 
only the difference between these objects but that the existence of one of them 
does not guarantee the existence of another were constructed and studied in \cite{BB, Mi,JT}. 
Curiously, despite that the authors of \cite{JT} claimed that they 
``give a complete description of relations between observable and natural measures'', 
a number of situations were not taken into account there. 
The most important among those omitted are the case when the limit 
measure $\t\mu^{\rm nat}$ is singular with respect to the reference measure and 
the case when $\t\mu^{\rm nat}$ does not exist. We shall study both these situations 
in the present paper.

It is worth noting that the so called operator approach (see e.g. \cite{Bl}) provides very 
effective tools for the analysis of the natural measures which makes them more preferable 
from the applied point of view than the observable ones.

\begin{definition} {\em Support} of a measure $\mu$, which is denoted by $\supp(\mu)$,  
is a union of all points $x\in X$ satisfying the property that each open neighborhood of 
the point $x$ has strictly positive $\mu$-measure. For a given natural measure $\t\mu$ 
we denote by $S$ its support and by $m_S$ the conditional probability measure constructed 
from the reference measure $m$ on the set $S$.
\end{definition} 

To simplify the presentation we assume that if a sequence of points $\{x_n\}$, which belongs 
either to the set $S$ or to its complement, converges to a point $x$ as $n\to\infty$, then the 
sequence of their images $\{Tx_n\}$ has a limit which in general may differ from $Tx$. 

The condition $m(S)>0$ implies the existence of the {\em conditional} measure 
$m_S$, namely $m_S(A):=m(A\cap S)/m(S)$, otherwise if $m(S)=0$ one needs 
additional assumptions for the existence of the conditional measure (see a general 
construction related to measurable partitions, e.g. in \cite{Ro}). 
Therefore to avoid these difficulties we assume that in the situations under consideration 
the measure $m_S$ is always well defined and that the map $T$ is nonsingular with 
respect to this conditional measure. 
In the examples which we shall consider this is indeed the case. 

Combining the results obtained in \cite{BB,Mi,JT,Hu} we get the following information  
about the connections between natural and observable measures.

\begin{theorem}\label{t:bb} \cite[theorem 2.1]{BB}, \cite[theorem 2.4]{JT}, \cite[theorem 1]{Hu}  
Let $T$ be a measurable map, having no wandering sets of positive $m$-measure,
and let there exist $\t\mu^{\rm obs}$ with the open basin of attraction $U$ of positive 
$m$-measure. Then there exists the natural measure $\t\mu^{\rm nat}$ with the same 
basin of attraction and $\t\mu^{\rm nat}=\t\mu^{\rm obs}$. Conversely, let there exist  
a natural $T$-invariant measure $\t\mu^{\rm nat}$ with the open basin of attraction $U$ 
of positive $m$-measure. Then $m_U \ll \t\mu^{\rm nat}$ implies that $\t\mu^{\rm nat}$ 
is ergodic and observable (i.e. $\t\mu^{\rm obs}=\t\mu^{\rm nat}$).
\end{theorem}

The importance of the nonexistence of wandering sets of positive $m$-measure to 
the construction of invariant measures was first observed in \cite{Hu} and later it has been 
shown that conditions of this sort are necessary for the presence of invariant measures 
absolutely continuous with respect to $m$-measure (see e.g. \cite{HK}). In our setup 
it is simpler to consider a slight generalization of this notion -- wandering measures 
absolutely continuous with respect to the reference measure. 

It turns out that the assumption of continuity of the map $T$ simplifies a lot the construction 
of the natural measure.

\begin{theorem}\label{t:mi} \cite[theorem 2.1]{Mi}  
Let $T$ be a continuous map and let the limit measure 
$\t\mu:=\lim_{n\to\infty}\frac1n\sum_{k=0}^{n-1}T_*^k m$ exist and be ergodic 
with respect to $T$. Then $\t\mu$ is the natural measure with the 
basin of attraction $\supp(m)$.
\end{theorem}

An alternative approach to is based on the idea of genericity. For example, 
the ergodicity assumption used in all ergodic theorems mentioned in 
Section~\ref{s:intro} may be justified by the well known Oxtoby-Ulam 
result \cite{OU}, according to which a generic volume preserving 
homeomorphism of a compact manifold is ergodic. This explains to some 
extent reasons why this is exactly the case considered in a vast majority of 
textbooks on dynamical systems theory. 
In turn ergodicity by the Birkhoff ergodic theorem means that 
almost all points (with respect to the volume measure) are typical in 
this case. In a recent paper \cite{AA} this reasoning was extended 
to generic continuous maps {\em without} the assumption of the 
volume preservation and it has been shown that for a generic map 
the Birkhoff average converges almost everywhere, but the limit value 
may depend sensitively on the initial point. This disproves 
a conjecture by D.~Ruelle \cite{Ru} who expected that generically 
those averages should diverge, which he called {\em historical 
behavior}. A different approach to this question together with a 
comprehensive review of corresponding results may be found in \cite{AP}, 
see also \cite{CM,GK,Ca}. 
Anyway, this discussion demonstrates that without some specific assumptions 
about the dynamics there is no hope to obtain positive results in this direction.

\section{Regular case}\label{s:inv+}

By the ``regular case'' we mean that the map $T$ is good enough and the natural measure 
$\t\mu$ is ``finer'' than the reference measure $m$. To be precise we make the following
{\em standing assumption} in this Section: the transfer-operator $T_*$ is continuous 
at the point $\t\mu$ and $S:=\supp(\t\mu)\subseteq\supp(m)$. Observe that this property 
does not require the map $T$ to be continuous.

\subsection{Main results}\label{s:results}

Throughout the paper we fix a compact metric phase space $(X,\rho)$ equipped with 
a Borel $\sigma$-algebra of measurable sets $\cB$ and a probabilistic 
reference measure $m$ on this $\sigma$-algebra. To have a simple picture in mind 
a reader may think that we deal with the unit interval $X:=[0,1]$ with the uniform metric, 
the standard Borel $\sigma$-algebra and the Lebesgue measure $m$ on it. 
Let us formulate our main results in the classical setting when invariant 
measures of a dynamical system under study are present. Namely, we 
assume that the natural measure $\t\mu$ exists. Additionally we assume that 
the map $T$ is nonsingular not only with respect to the reference measure $m$, 
but with respect to the conditional reference measure $m_S$. 

\begin{theorem} (Necessary and sufficient conditions) \label{t:nsc}
The property $m(Z_{\t\mu})\cdot m_S(Z_{\t\mu})=1$ is equivalent to 
the following 3 assumptions %
\begin{itemize}
\item [(i)] $\frac1n\sum_{k=0}^{n-1}T_*^k\mu\toas{n\to\infty}\t\mu \quad \forall 
    \mu\in\cM(m)\cup\cM({m_S})$, 
\item [(ii)] the limit measure $\t\mu$ is ergodic, 
\item [(iii)] there are no wandering measures in $\cM(m_S)$.
\end{itemize}
\end{theorem}

This result is in fact stronger than what we have discussed in the Introduction, 
since we claim additionally that the conditional reference measure $m_S(Z_{\t\mu})=1$. 

A somewhat unusual assumption (iii) of Theorem~\ref{t:nsc} may be 
replaced by another one looking much simpler, but in fact being stronger.

\begin{remark} \label{r:iii} The assumption (iii) holds if \par
{\em (iii')} $\t\mu\in\cM(m_S)$.
\end{remark}

Indeed, if (iii') holds true then the presence of an wandering measure absolutely 
continuous with respect to the measure $m_S$ would contradict to the definition 
of the invariant measure $\t\mu$. 
This new assumption is much stronger and does not allow to study situations 
with singular invariant measures, but it is easier to control, which we shall 
use in Section~\ref{s:self} dedicated to self-consistent dynamical systems. 
Note that the last claim made in Theorem~\ref{t:bb} is very similar to this 
assumption. 

\bigskip

In the the case of continuous maps above conditions may be considerably simplified.

\begin{theorem} (Case of continuous maps) \label{t:suf-c} Let $T:X\to X$ 
be a continuous map. Then the $m$-a.e. convergence in (\ref{e:trj-conv}) 
is equivalent to the following 3 conditions: %
\begin{itemize}
\item [(i)] $\frac1n\sum_{k=0}^{n-1}T_*^k m\toas{n\to\infty}\t\mu$,
\item [(ii)] the limit measure $\t\mu$ is ergodic, 
\item [(iii)] there are no wandering measures in $\cM(m_S)$.
\end{itemize}
Moreover, if $\t\mu\ll m$, then the conditions (ii), (iii) may be omitted altogether.
\end{theorem}

\subsection{Analysis of the examples}\label{s:example}

The main aim of this Section is to discuss in detail ergodic properties of 
the dynamical systems introduced as examples in the Introduction and 
their simple modifications.

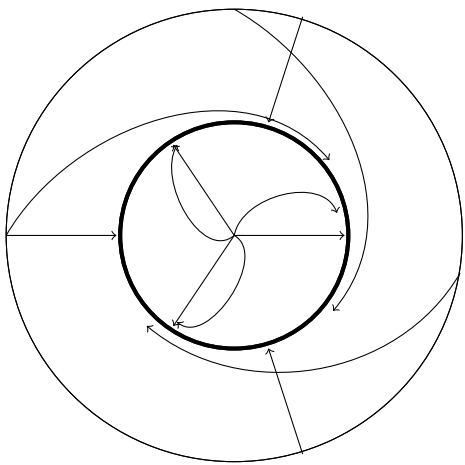
\begin{figure}
\begin{center}\begin{tikzpicture}
\put(-120,0){\draw circle (3cm); \draw[line width=1.5pt] circle (1.5cm);
               \draw [->]  (-3, 0) to [out=60, in=130] (1.25, 1.0);
               \draw [->]   (0, 3) to [out=-30, in=50] (1.3, -1);
               \draw [->]  (2.97, -0.5) to [out=-120, in=-40] (-1.15, -1.2);
               \draw [->]  (0, 0)  to [out=-30, in=-40] (-0.75, -1.15);
               \draw [->]  (0, 0)  to [out=-140, in=-120] (-0.75, 1.2);
               \draw [->]  (0, 0)  to [out=80, in=110] (1.35, 0.3);
               } 
\put(120,0){\draw circle (3cm); \draw[line width=1.5pt] circle (1.5cm);
               \draw [->]  (-3, 0) -- (-1.55,0); \draw [->]  (0.9, 2.9) -- (0.45,1.5);
               \draw [->]  (0.9,-2.9) -- (0.45,-1.5);
               \draw [->]  (0,0) -- (1.45,0);
               \draw [->]  (0,0) -- (-0.80,1.2); \draw [->]  (0,0) -- (-0.80,-1.2);
               } 
\end{tikzpicture} \end{center}
\caption{Phase portraits of the dynamical systems in the examples \ref{ex:1} (left) 
and \ref{ex:3} (right)}\label{f:examples}
\end{figure}

Since the situation with the one-dimensional contracting map $Tx:=x/2$ 
(mentioned in the Introduction) is trivial, we start with the example~\ref{ex:1} 
(the corresponding phase portrait is sketched in Figure~\ref{f:examples} (left)). 
Recall that in this example we consider a family of maps from the unit disc 
$X:=\{(\phi,R):~0\le\phi<2\pi,~0\le R\le1\}$ into itself defined in the polar 
coordinates $(\phi,R)$ by the relation: %
$$ {T(\phi,R):=(\phi+2\pi\alpha+\beta(R-r)~{\rm mod~2\pi},~\gamma(R-r)+r)} .$$ %
This simple example is very instructive since choosing a different quadruple 
of admissible parameters $(\alpha,\beta,r,\gamma)\in(0,1)^4$ or making a slight 
modification of this system we are able 
to construct illustrations and counter-examples to our main results formulated 
in the previous Section. 

Let the reference measure $m$ be chosen as the 2-dimensional Lebesgue measure 
on $X$ normalized by $\pi$, $C:=\{R=r\}$ -- the centered circle of radius $r$, 
and let $m_C$ be the 1-dimensional Lebesgue measure on $C$, normalized by $2\pi r$, 
i.e. the normalized conditional 2-dimensional Lebesgue measure on the 1-dimensional set $C$. 

Obviously for each admissible quadruple $(\alpha,\beta,r,\gamma)$ the system 
$T,X$ possesses a single attractive set $C$ (with respect to the Euclidean metric 
$\rho(\cdot,\cdot)$ on $X$) and the measure $m_C$ is $T$-invariant. 

\begin{proposition}\label{p:ex1} \begin{itemize}
\item [(a)] $\rho(T^nx,C)\toas{n\to\infty}0$ for all $x\in X \setminus \{0\}$.
\item [(b)] $(T,C,m_C)$ is ergodic if and only if $\alpha$ is irrational. 
\item [(c)] $m(Z_{m_C})=1$ whenever $\alpha$ is irrational. 
\item [(d)] The measure $m_C$ is natural for each admissible $\alpha$.
\end{itemize} \end{proposition}

All these claims are more or less straightforward and we leave the proof for the reader.

Observe that by (d) even for rational values of the parameter $\alpha$ each 
$m$-smooth measure converges to $m_C$ in Cesaro means. Moreover, 
making an arbitrary $C^\infty$-small perturbation to the maps from 
this family (slowing down the rotation around $S$ in its small neighborhood) one gets 
a system, for which $(T,C,m_C)$ is non-ergodic, i.e. $\alpha$ 
is rational, but the property $m(Z_{m_C})=1$ remains valid. This explains, why the 
properties $m(Z_{m_C})=1$ and $m_C(Z_{m_C})=1$ should be considered 
separately, despite that they look very similar.

Let us show that eliminating the rotation outside of the attractor $S$ in the example~\ref{ex:1} 
we come to the situation when $m(Z_{m_C})=0$ while $m_C(Z_{m_C})=1$. 

\begin{example}\label{ex:3}
$$ T(\phi,R) := \left(\function{\phi + 2\pi\alpha & \mbox{if } R=r,  \\
                                             \phi & \mbox{if } R\ne r}, ~~\gamma(R-r) \right)  $$
\end{example}

The corresponding phase portrait is sketched in Figure~\ref{f:examples} (right). 
It is easy to see that here we have the uniform global convergence to the support 
of the measure $m_C$ (albeit without the property (i)), but for any point 
$x:=(\phi,R)$ with $R\ne r$ the ergodic average $\frac1n \sum_{k:=0}^{n-1}f(T^nx)$
converges to $f(\phi,r)$ rather than to $\int f dm_C$.  

%The map in the example~\ref{ex:0} is continuous and thus obviously satisfies our 
%regularity assumption. On the other hand, 

\subsection{Proofs}\label{s:proofs}
We start with the connections between the point-wise convergence described 
by the relation~(\ref{e:star}) and the convergence of images of measures 
under the action of the dynamical system, which is formulated as the property 
(i) in above theorems. 
Observe that the property (\ref{e:star}) is equivalent to the existence (and 
hence uniqueness) of the observable measure $\t\mu^{\rm obs}$ having  
the entire space $X$ as the basin of attraction.
 
Fix a measurable map $T$ from a compact measurable metric space $(X,\cB)$ into itself. 
 
\begin{definition} A point $x\in X$ is called {\em typical} for a probability measure $\mu$ if 
$$  \frac1n \sum_{k=0}^{n-1}T_*^k\delta_x \toas{i\to\infty}\mu $$
in the weak sense. 
Recall that the set of typical points for the measure $\mu\in\cM$ is denoted by $Z_\mu$.
\end{definition}

In these terms an observable measure is characterized by the property that there 
is an open set of typical points of positive $m$-measure. 

\begin{lemma}\label{l:impl-inv} Let $\mu$ be a limit point for the sequence 
of measures $\frac1n \sum_{k=0}^{n-1}T_*^k\delta_x$ for some point 
$x\in X$ and let the map $T_*:\cM(X)\to\cM(X)$ be continuous at the 
measure $\mu$. Then $T_*\mu=\mu$. 
\end{lemma}
\proof By the assumption of Lemma there exists a sequence of indices 
$\{n_i\}_i\toas{i\to\infty}\infty$ such that 
$$ \mu_{n_i,x}:= \frac1{n_i} \sum_{k=0}^{n_i-1}T_*^k\delta_x \toas{i\to\infty}\mu .$$
On the other hand, for each continuous function $\phi$ we have 
$$ |(\delta_x - T_*\delta_x)(\phi)| \le 2|\phi|_\infty .$$
Since the map $T_*$ is continuous at the measure $\mu$ one can interchange 
the action of the transfer-operator $T_*$ and the operation of passing to the limit. 
Therefore  %
\bea{ T_*\mu \a= \lim_{i\to\infty} T_*\mu_{n_i,x} \\
                 \a= \lim_{i\to\infty} (\mu_{n_i,x} - \frac1{n_i}(\delta_x - T_*^{n_i}\delta_x)) \\
                 \a= \lim_{i\to\infty} \mu_{n_i,x} = \mu .} %
Observe that without the continuity assumption one cannot make a conclusion about 
the invariance of the limit measure. \qed

One of the key points of our argument is the following result describing metric 
properties of basins of attraction of general non necessarily invariant measures. 

\begin{lemma}\label{l:mes-typ} If $\mu$ is an ergodic $T$-invariant measure, then 
$\mu(Z_\mu)=1$. Otherwise, if $\mu$ is non-ergodic or ergodic but non-invariant 
and the map $T_*:\cM(X)\to\cM(X)$ is continuous at the measure $\mu$, then $\mu(Z_\mu)=0$.
\end{lemma}
\proof If the measure $\mu$ is an ergodic $T$-invariant measure this claim is 
a trivial consequence of the Birkhoff ergodic theorem. Indeed, it is enough to 
apply ergodic theorem for the test-function defined as the indicator function 
of the set $Z_\mu$. 

It remains to prove that in the opposite situation the set of typical points 
is of zero $\mu$-measure. By Lemma~\ref{l:impl-inv} a measure having even a 
single typical point at which the transfer-operator is continuous inevitably 
has to be an invariant measure. Therefore by the Ergodic Decomposition Theorem 
(see e.g. \cite{Si}) there is a probability measure $\eta$ on the set $\cM_T^{\rm erg}$ of all 
ergodic $T$-invariant measures $\nu$, such that for each measurable set $B\in\cB$ 
we have %
\beq{e:erg-decomp}{ \mu(B) = \int \nu(B) ~d\eta(\nu) .} %
By the first part of the proof, $\nu(Z_\nu)=1$ for each ergodic measure $\nu$. 
On the other hand, by the definition of the set of typical points  
$Z_\nu\cap Z_{\nu'}=\emptyset$ 
for any two measures $\nu,\nu'$, which implies that 
$$ Z_\mu \subseteq X \setminus \cup_{\nu\in\cM_T^{\rm erg}}Z_\nu .$$
Applying (\ref{e:erg-decomp}), we get 
$$ \mu(\cup_{\nu\in\cM_T^{\rm erg}}Z_\nu)=1 .$$
Therefore 
$$ \mu(Z_\mu) \le \mu( X \setminus\cup_{\nu\in\cM_T^{\rm erg}}Z_\nu) = 0 ,$$
which finishes the proof. \qed

\begin{remark} If $T\in C^0(X,X)$ then $T_*$ is continuous at any measure $\mu$, 
but otherwise $T_*$ might be discontinuous even at the invariant measure.
\end{remark}

Indeed, consider  
\begin{example}
$X:=[0,1]$ and 
$Tx:=\function{x/2+1/4&\mbox{if } 0\le x\le1/2 \\  2x-1 &\mbox{otherwise}}$. 
\end{example}
Then the transfer-operator $T_*$ is discontinuous at the $T$-invariant measure $\delta_{\frac12}$. 

\begin{lemma}\label{l:path-mes} The property $m(Z_{\t\mu})=1$ implies (i).
\end{lemma}
\proof By definition for a
$m$-a.a. point $y\in X$ the weak convergence of the sequence of measures 
$\frac1n\sum_{k=0}^{n-1}\delta_{T^ky}\to\t\mu^{\rm obs}$
takes place, i.e. for any continuous function $\phi:X\to\IR^1$ and $m$-a.a. $y\in X$ 
we, using the Lebesgue dominated Convergence Theorem, get 
$$ \int \phi(x)~d\left(\frac1n\sum_{k=0}^{n-1}\delta_{T^ky}\right)
  = \frac1n\sum_{k=0}^{n-1}\phi(T^ky) \toas{n\to\infty}
    \int \phi~d\t\mu^{\rm obs} .$$
Choose an absolutely continuous (with respect to $m$) measure $\mu\in\cM(X)$ 
and consider Cesaro averages of its images: 
$\mu_n:=\frac1n\sum_{k=1}^{n-1}{T^*}^k\mu$. 
Then using the above convergence and the absolute continuity of the
measure $\mu$ we get %
\bea{ \int \phi~d\mu_n
  \a = \int \phi ~d\left(\frac1n\sum_{k=0}^{n-1}{T^*}^k\mu\right) \\
  \a = \int \frac1n\sum_{k=0}^{n-1}\phi(T^kx) ~d\mu \\
  \a \toas{n\to\infty} \int \left(\int \phi~d\t\mu^{\rm obs}\right)~d\mu
     = \int \phi~d\t\mu^{\rm obs} ,}% 
which proves the assertion. \qed

Consider now the questions about the necessity of ergodicity and the absence of wandering sets. 

\begin{lemma}\label{l:erg-nec} The property $m(Z_{\t\mu})=1$ implies ergodicity, i.e. (ii).
\end{lemma}
\proof By Lemma~\ref{l:mes-typ} a non-ergodic invariant measure cannot have 
a full $m$-measure set of typical points, which implies the claim.\qed

Despite of this result even the uniqueness of the measure $\t\mu^{\rm obs}$ 
does not imply that it is ergodic. Indeed, consider the
following example (see \cite{BB} for details):%
\beq{ex:GiGi}{
   T x := \function{(1-\sin(\pi x-\pi/2))/2 & \mbox{if } 0< x < 1,  \\
                                             x & \mbox{if } x\in\{0,1\} } .}%
One can easily show that the locally maximal attractor in this
example consists of two fixed points at 0 and 1 and that
$\frac1n\sum_{k=1}^{n-1}\delta_{T^n x}
   \toas{n\to\infty}\frac12(\delta_0+\delta_1) = \t\mu^{\rm obs}$
for any initial point $x\in(0,1)$. 
On the other hand, this measure is nonergodic, since the points
0 and 1 are fixed points. 

According to Lemma~\ref{l:mes-typ} for the analysis of the regular case 
we need to study only properties of ergodic invariant measures with a nontrivial 
support $S$. Thus the set $S$ is forward invariant.

\begin{lemma}\label{l:wand-nec} Let $\t\mu$ be an ergodic invariant measure.  
Then the property $m_S(Z_{\t\mu})=1$ implies the absence of wandering 
measures in $\cM(m_S)$.
\end{lemma}
\proof Assume from the contrary that there exists an wandering measure 
$\mu\in\cM(m_S)$. Denote $\t{S}:=Z_{\t\mu}\cap S$. 
Then by the assumption the measure 
$\mu_{\t{S}}$ also belongs to the set $\cM(m_S)$. Hence 
$$ \frac1n\sum_{k=0}^{n-1}T_*^k \mu_{\t{S}} \toas{n\to\infty}\t\mu .$$
By the wandering property the measures $T_*^k \mu_{\t{S}}$ are mutually 
singular for different $k$, which together with the forward invariance of $S$ 
contradicts to the fact that the support of the limit measure $\t\mu$ coincides 
with the set $S$. \qed

It is worth noting that this result is very similar to the analysis of transformations 
without wandering sets of positive measure in \cite{Hu}(\S 10). %Section 

\bigskip

\n{\bf Proof of Theorem~\ref{t:nsc}}. Collecting together the results of lemmas 
\ref{l:path-mes}, \ref{l:erg-nec}, \ref{l:wand-nec} 
we get the first part of the Theorem related to the necessity of the 
assumptions (i)--(iii). 

Now we turn to the second part of the proof of Theorem~\ref{t:nsc} 
and consider in detail the set 
of points converging under dynamics to the support $S$ of the measure $\t\mu$: 
$$ Y:=\{x\in X:~~\lim_{t\to\infty} \rho(T^tx,S)=0\} .$$

\begin{lemma}\label{l:mes-trace} Let the condition \ref{t:nsc}(i) hold true. Then $m(Y)=1$.
\end{lemma}
\proof Assume from the contrary that this claim does not hold. Then there exists 
a subset %
\beq{e:bad-set}{A:=\{x\in X:~~\liminf_{t\to\infty} \rho(T^tx,S)>0\}.} % 
of positive $m$-measure. Denote by $m_A$ the conditional measure induced by 
the reference measure $m$ on the set $A$. 
Since $m(A)>0$ this measure is absolutely continuous with 
respect to $m$ and hence by \ref{t:nsc}(i) we have 
$$ \frac1n\sum_{k=0}^{n-1}T_*^k m_A\toas{n\to\infty}\t\mu ,$$
which contradicts to the definition of the set $A$. \qed

\begin{definition} We say that a point $y_x\in S$ is {\em weakly tracing} a point $x\in X$ 
if the following limits (in the weak sense) exist and coincide 
$$  \lim_{n\to\infty}\frac1n \sum_{k=0}^{n-1}T_*^k\delta_x  = 
     \lim_{n\to\infty}\frac1n \sum_{k=0}^{n-1}T_*^k\delta_{y_x} .$$
\end{definition}

In other words, trajectories starting from the points $x\in X$ and $y_x\in S$ 
have the same statistics. 
Obviously the map $x\to y_x$ is not injective, however, applying the same argument 
as in the proof of Lemma~\ref{l:mes-trace}, we deduce that 

\begin{corollary}\label{c:tracing} For $m$-a.e. point $x\in Y$ there is a weakly 
tracing point $y_x\in S$.
\end{corollary}

Observe that one cannot use here the stronger point-wise tracing property 
$$ \rho(T^nx,T^ny_x) \toas{n\to\infty}0 $$
Indeed, consider a contractive system whose trajectories are wrapping around 
the limit circle. The dynamics on the circle is supposed to be a pure irrational rotation 
like in the Example~\ref{ex:1}. Then if the rate of convergence to the attractor is 
slow enough, e.g. of order $1/n$, then the outer points do not possess the 
point-wise tracing counterparts on the attractor. 

\begin{definition} Denote by  $\tilde{Y}$ the subset of points from the set $Y$ 
weakly traced by $\t\mu$-``typical'' trajectories on the support of the measure 
$\t\mu^{\rm nat}$. \end{definition} 

Our aim is to check that this set of full $m$-measure. To this end one is tempted 
to extend the argument used in the previous proof along the following lines.

Assume from the contrary that there exists a subset $B\subseteq Y$ of positive $m$-measure 
such that the trajectories starting from this set are traced by non-typical points from $S$. 
By the ergodicity assumption \ref{t:nsc}(ii) according to Lemma~\ref{l:mes-typ} we have $\t\mu(Z_{\t\mu})=1$. Hence the 
$\t\mu$-measure of the set of non-typical points $S\setminus Z_{\t\mu}$ is zero. 
Denote by $m_B$ the conditional measure constructed from $m$ on the set $B$. 
Since $m(B)>0$ this measure is absolutely continuous with respect to $m$ and hence by 
\ref{t:nsc}(i) we have 
$$ \frac1n\sum_{k=0}^{n-1}T_*^k m_B\toas{n\to\infty}\t\mu .$$
On the other hand, since by construction all points from the set $B$ are traced by 
$\t\mu$-non-typical points, the limit measure should be supported by the set 
$S\setminus Z_{\t\mu}$, which seems to be a contradiction. 

Unfortunately the set of non-typical points being $T$-invariant needs not to be 
a compact set. Therefore the limit measure may be supported by a larger set of 
positive $\mu$-measure. This explains that some additional assumptions are 
necessary in order to overcome this difficulty. 

One of the possibilities is to assume that the measures $m_S$ and $\t\mu$ 
are equivalent (see Remark~\ref{r:iii}). However we prefer to use a less restrictive 
assumption about the absence of wandering measures. 

\begin{lemma}\label{l:tracing} Let the conditions \ref{t:nsc}(i,ii,iii) hold true and 
let $m(S)>0$. Then $m(\tilde{Y})=1$.
\end{lemma}
\proof Denote $\tilde{Z}:=Z_{\t\mu}\cap S$. From the previous results we have 
$\t\mu(\tilde{Z})=1$, but this fact alone does not contradict to $m(\tilde{Z})<1$. 
If the latter inequality takes place, then the set $B:=S\setminus \tilde{Z}$ is of 
positive $m$-measure. By definition only $\t\mu$-non-typical points belong to the set $B$. 
Then this together with the assumption (i) and the ergodicity of $\t\mu$ imply that 
the measure $m_B$ is wandering, which contradicts to the assumption (iii). 

Consider now what happens outside of the set $S$. Denote $Z:=Z_{\t\mu}\setminus S$. 
Assume from the contrary that there exists a subset $C\subseteq(X\setminus S)\cap Y$ 
of positive $m$-measure consisting only of $\t\mu$-non-typical points. Since $m(Y)=1$ 
by Lemma~\ref{l:mes-trace}, the intersection with $Y$ does not change the measure of 
the set $C$. By Corollary~\ref{c:tracing} for each point $x\in C$ there is a tracing 
point $y_x\in S$. Decompose the set $C$ into two parts: $C_0$ consisting of points 
starting from which a trajectory hits the set $S$ in a finite number of steps, and $C_1$ 
for which trajectories only converge to $S$. By the 1-st part of the proof $m(C_0)=0$ 
(as a union of pre-images of the set of $m$-measure zero), while to prove that $m(C_1)=0$ 
one uses additionally the property (i).  \qed

\begin{lemma}\label{l:tracing-zero} Let the conditions \ref{t:nsc}(i,ii,iii) hold true,   
and let $m(S)=0$, but $m_S$ be well defined. Then $m(\tilde{Y})=1$.
\end{lemma}
\proof The only difference with the proof of Lemma~\ref{l:tracing} is that one needs 
to use additionally the assumption that the property~(\ref{t:nsc}(i)) holds not only 
with respect to the reference measure $m$, but also with respect to the conditional 
measure $m_S$, which is singular with respect to $m$ in this case. Therefore we 
omit details.  \qed 

\n{\bf Final part of the proof of Theorem~\ref{t:nsc}}. Collecting the results obtained in 
lemmas~\ref{l:mes-trace}--\ref{l:tracing-zero} we get that $m(Z_{\t\mu})=1$. 
To obtain the claim that $m_S(Z_{\t\mu})=1$ one needs to consider the 
restriction of the dynamical system to the forward invariant set $S$ 
(similarly to the proof of the necessary conditions). \qed 

It remains to {\bf prove Theorem~\ref{t:suf-c}}. Here the item (i) follows from 
Theorem~\ref{t:mi}. The consideration of items (ii) and (iii) is exactly the same 
as in the proof of the previous result. \qed

\section{Irregular case}\label{s:irr}
This section is dedicated to the situation when the transfer-operator $T_*$ is discontinuous 
at the limit measure $\t\mu:=\lim_{n\to\infty} \frac1n\sum_{k=0}^{n-1}T_*^k\mu$. 
Normally this happens when the dynamical system has no invariant measure. 
A typical example of this sort may be obtained by the following simple modification 
of the example~\ref{e:circle}:% 
\begin{example}\label{ex:4}
\beq{e:circle-mod}{T(\phi,R):=\function{
     (\phi+2\pi\alpha+\beta(R-r)~{\rm mod~2\pi},~\gamma(R-r)+r)  &\mbox{if } r(R-r)\ne0 \\
     (\phi+2\pi\alpha~{\rm mod~2\pi},~(1+r)/2) &\mbox{otherwise}.}
}  %
\end{example}
This example and its modifications are very instructive, but to have a simpler 
picture in mind let us consider the following one-dimensional map:

\begin{example} \label{ex:2}
$X:=[0,1)$ and $Tx:=\function{1-c &\mbox{if } x=0 \\ x^2 &\mbox{otherwise}}$.  
\end{example}
If $0<c<1$, then this map has no invariant measures, but each absolutely continuous 
measure converges to $\delta_0$; otherwise if $c=0$ the situation changes 
drastically: the limit measure is invariant. 

One can show that in the examples above for any smooth probability measure $\mu$ 
the sequence of measures $\frac1n\sum_{k=0}^{n-1}T_*^k\mu$ converges weakly 
to a certain limit measure $\t\mu$ on the circle $\{R=r\}$ or at the origin, 
but this measure is no longer invariant (if $c\ne0$ in the example~\ref{ex:2}). 
Instead of the convergence to the limit circle or to a fixed point one can consider 
situations when the only attractor of the dynamical system $(T,X)$ is a Cantor set $K$. 
Then modifying the map $T$ only on the set $K$ in the same way as above 
(e.g. $\t{T}|_{X\setminus K}\equiv T|_{X\setminus K},~~ \t{T}K\ne K$) 
one gets nontrivial examples of the discontinuity of the transfer-operator $\t{T}_*$ 
at the limit measure.

In order to apply arguments similar to those elaborated in Section~\ref{s:inv+} 
we need to clarify and extend the notion of the invariant measure. 
To this end we make use of the equivalence of the ergodicity and the fulfillment 
of the Birkhoff Theorem (see e.g. \cite{Si}). 

\begin{definition} Motivated by the Birkhoff ergodic theorem, we say that a measure 
$\mu\in\cM$ is {\em weakly ergodic} if $\mu(Z_\mu)=1$. 
\end{definition}

Indeed, this property obviously holds in the case of a conventional ergodic invariant 
measure. To demonstrate the reason for this extension, consider the the example~\ref{ex:2}. 
When $0<c<1$ this map has no invariant measures, 
but the $\delta$-measure at the origin is weakly ergodic. More sophisticated 
situations, concerning the so called nonlinear Markov chains, will be discussed 
in Section \ref{s:self}.

The assumption about the discontinuity of the transfer-operator $T_*$ at the 
measure $\t\mu$ implies obviously the discontinuity of the map $T$. 
Nevertheless making some additional assumptions about local properties of 
the map $T$ in the neighborhood of the support of the measure $\t\mu$ 
in principle one can get an analogue of the necessary part of Theorem~\ref{t:nsc}. 
However those assumptions look somewhat clumsy and therefore we restrict 
ourselves only to sufficient conditions.

\begin{theorem} \label{t:suf-irr} The assumptions
\begin{itemize}
\item [(i)] $\frac1n\sum_{k=0}^{n-1}T_*^k\mu\toas{n\to\infty}\t\mu \quad \forall 
    \mu\in\cM(m)\cup\cM({m_S})$, 
\item [(ii)] the limit measure $\t\mu$ is weakly ergodic, 
\item [(iii)] there are no wandering measures in $\cM(m_S)$,
\end{itemize}
imply that $m(Z_{\t\mu})\cdot m_S(Z_{\t\mu})=1$.
\end{theorem} 

The proof follows basically the same scheme as in the regular case with the necessary 
usage of the weak ergodicity property instead of references to the ergodic theorem, e.g.: 

\begin{lemma} %\label{l:mes-typ-irr} 
If $\mu$ is weakly ergodic $T$-invariant measure, then $\mu(Z_\mu)=1$.
\end{lemma}
\proof The claim is an immediate consequence of the definition of the weakly ergodic measure. \qed

This result is an extension of Lemma~\ref{l:mes-typ}. Further on the statements of 
Lemmas~\ref{l:mes-trace} -- \ref{l:tracing-zero} made in the regular case 
do not make use of the continuity of the map $T$ and thus remain valid in the discontinuous 
setting. Collecting the corresponding results we get the proof of Theorem~\ref{t:suf-irr}. \qed

\section{Self-consistent dynamical systems}\label{s:self}

In this section we discuss a more complicated and not well studied case when 
the dynamics depends not only on the current point in the space $X$ but on the 
current statistics of the system as well. Let $\{T_\mu\}_{\mu\in\cM(X)}$ be a family 
of maps from a compact measurable space $X$ into itself, parametrized 
by probability measures $\mu\in\cM(X)$.

\begin{definition} By a {\em self-consistent dynamical system} we mean 
a skew product map $\cT(x,\mu):=(T_\mu x,~ (T_\mu)_*\mu)$ 
%$\left(\matrix{x\to T_\mu x\cr \mu\to (T_\mu)_*\mu \cr}\right)$  
acting in the direct product space $X\times\cM(X)$. 
\end{definition}

The idea of deterministic self-consistent dynamical systems was introduced 
by K.~Kaneko \cite{Ka} in order to approximate the dynamics of large systems 
in terms of a  mean-field type perturbation of an isolated sub-system.  
Despite a large number of attempts to study such systems there are only a few 
situations when a complete mathematical treatment was successful. See \cite{Ke,BKZ} 
and further references to known numerical results therein. Discussion of similar 
questions in true random setting may be found e.g. in \cite{Bu}. 

Let us give a couple of seemingly trivial examples. Let $X:=[0,1]$ and denote by 
$E_\mu:=\int x d\mu$ the mathematical expectation over a probabilistic measure 
$\mu$ on $X$. 

\begin{example} \label{ex:sc-1} (Additive perturbation) 
$T_{\ep,\mu}x:=Tx + \ep E_\mu$~mod 1, where $\ep\in[0,1]$. 
\end{example}

\begin{example} \label{ex:sc-2} (Multiplicative perturbation) 
(a) $T_\mu x:=x\cdot E_\mu$~mod 1, \\ (b) $T_\mu x:=x/ E_\mu$~mod 1 
(here we assume that $1/0$ mod 1 = 0).  
\end{example}

The example~\ref{ex:sc-1} with $Tx:=1-2|x-1/2|$ was studied in \cite{Ke}, where it has been shown that 
images of any probability measure absolutely continuous with respect to the Lebesgue 
measure $m$ on $X$ converge in Cesaro means to $m$ provided the parameter 
$\ep$ be small enough. This is one of a very few cases where the stability of the 
SRB measure with respect to mean-field type perturbations has been proven. 

The example~\ref{ex:sc-2} represents the situation of multiplicative perturbations 
which was not studied earlier.  
In distinction to the case (a), where the dynamics is trivial: there are only two invariant 
measures stable one $\delta_0$ and unstable $\delta_1$, the dynamics in 
the case (b) is much more complicated: there are infinitely many mutually 
singular probabilistic invariant measures, including the Lebesgue measure. Here the 
invariance means that the measure is preserved under dynamics (this explanation is 
necessary since the system itself is defined in the product space $X\times\cM(X)$).

In the present setting, even under the assumption that all maps $T_\mu$ are continuous, 
one cannot apply the ergodic theorem directly, since at each time step a 
different map is chosen. Additionally the conventional definition of ergodicity 
does no make much sense here. To overcome this difficulty we adapt the 
notion of the weak ergodicity to the setting of self-consistent dynamical systems.

\begin{definition} A measure $\mu$ is said to be {\em weakly ergodic} for the 
self-consistent system if the property %
\beq{e:trj-conv-sc}{\frac1n\sum_{k=0}^{n-1}f(\cT^k(x,\mu)) \toas{n\to\infty} \int f~d\mu }% 
holds for each continuous test function $f(\cdot,\cdot)$, depending only on the first 
argument $x$, and for $\mu$-a.e. $x\in X$. 
\end{definition}

Note that the Lebesgue measure $m$ is weakly ergodic in the example~\ref{ex:sc-2}(b), 
however it is not clear whether smooth probabilistic measures converge to $m$ even 
in a certain weak sense under the action of this self-consistent dynamical system. 

Let the transfer-operator $(T_{\tilde\mu})_*$ be continuous 
at the limit measure $\t\mu$ (see the assumption (i) below) and which is finer 
than the given reference measure $m$, i.e. $S:=\supp(\tilde\mu)\subseteq\supp(m)$. 
The definition of $\t\mu$-typical points also needs to be modified due to 
the more complex structure of the phase space ($X\times\cM(X)$ instead of $X$):
%\beq{e:star-sc}
$$ { Z_{\tilde\mu}:=\{x\in X:~~
       \frac1n\sum_{k=0}^{n-1}f(\cT^k(x,\mu)) \toas{n\to\infty} \int f~d\t\mu 
       \qquad \forall\mu\in\cM(m),~ \forall f\in C^0(X)  \} .}$$% 

The following result gives sufficient conditions for the application of the 
ergodic averaging in the entire space. Since in the case under consideration 
the point-wise dynamics depends sensitively on the choice of the initial 
measure, the construction of necessary conditions is not clear at the moment. 

\begin{theorem}\label{t:self-con} Let 
\begin{itemize}
\item [(i)] $\frac1n\sum_{k=0}^{n-1}\cT_*^k\mu\toas{n\to\infty}\tilde\mu \quad \forall 
    \mu\in\cM(m)\cup\cM({m_S})$, 
\item [(ii)] the limit measure $\tilde\mu$ be weakly ergodic,
\item [(iii')] $\tilde\mu\in\cM(m_S)$. 
\end{itemize}
Then $m(Z_{\tilde\mu})\cdot m_S(Z_{\tilde\mu})=1$.
\end{theorem}

Again as in the previous Section the usage of the notion of weakly ergodic measures 
allows to study the dynamics outside of the set $S:=\supp(\t\mu)$ and to prove the 
convergence to $S$ for $m$-a.a. initial points outside of $S$, while the convergence 
of statistics of trajectories starting from the points inside of $S$ follows from the 
assumption (iii').  \qed

As we already discussed in Section~\ref{t:nsc} (see Remark~\ref{r:iii}) 
the assumption (iii') is somewhat restrictive and one would prefer a weaker assumption: 

{\it (iii)} there are no wandering measures $\mu\in\cM(m_S).$

Unfortunately there are two obstacles here. First, at each time step one applies a 
different map from the family $T_\mu$ and thus the notion of the wandering measure 
needs to be modified.

\begin{definition} A probability measure $\mu$ is said to be {\em wandering} for the 
self-consistent system $\cT$ if the measures $\cT_*^n\mu$ and $\cT_*^k\mu$ are 
mutually singular for all $n\ne k\in\IZ^+$. 
\end{definition} 

A much more delicate point is that due to the same reason (different maps are applied 
at different time steps) the support of the limit measure $\tilde\mu$ needs not to be 
an invariant set for the maps from the family $T_\mu$. Therefore the situation with 
singular with respect to the conditional measure $m_S$ limit measures is out of control 
under the present approach. Nevertheless we expect that the assumption (iii) should 
be sufficient here. 

\section*{Acknowledgments}
The author is grateful to anonymous referees for helpful comments and 
suggestions which improved the quality of this paper.

%\newpage \small%\footnotesize


\begin{thebibliography}{99} 

\bibitem{AA} F. Abdenur, M. Andersson. 
  Ergodic theory of generic continuous maps, 
  Commun. Math. Phys. 318 (2013), 831-855. %arXiv:1201.0632v2 [math.DS]

\bibitem{AP} V. Araujo, V. Pinheiro. 
   Abundance of wild historic behavior, ergodic decomposition and generalized physical measures, 
   arXiv:1609.05356v1 [math.DS] %(40pp, 17 Sep 2016)

\bibitem{BKZ} J.-B. Bardet, G. Keller, R. Zweimuller. 
    Stochastically stable globally coupled maps with bistable thermodynamic limit,
    Comm. in Mathematical Physics 292 (2009), 237-270.

\bibitem{Bl} M. Blank. 
   Discreteness and continuity in problems of chaotic dynamics, Monograph, Amer. Math. Soc., 1997.

\bibitem{BB} M.  Blank, L. Bunimovich.  
   Multicomponent dynamical systems: SRB measures and phase transitions, 
    Nonlinearity, 16:1(2003), 387-401.
   % (math.DS/0202200) 

\bibitem{Bu} O. A. Butkovsky. 
      On Ergodic Properties of Nonlinear Markov Chains and Stochastic McKean--Vlasov Equations. 
      Theory Probab. Appl., 58:4, 661–674.

\bibitem{Ca} E. Catsigeras.  
       Ergodic Theorems with Respect to Lebesgue, 
       Trans. on Math. of the World Sci Eng Acad Soc, 10:12 (2012), 463 -479.

\bibitem{CM} J. Chaika, H. Masur. 
  There exists an interval exchange with a non-ergodic generic measure, 
  J. of Modern Dynamics (JMD), 9 (2015), 289-304. 

\bibitem{GK} K. Gelfert, D. Kwietnia. 
  On density of ergodic measures and generic points, 
  arXiv:1404.0456v2 [math.DS] %(32pp, 25 Aug 2015) 

\bibitem{HK} A.B. Hajian and S. Kakutani
  Weakly Wandering Sets and Invariant Measures,
  Trans. of AMS, 110:1(1964), 136-151.

\bibitem{Hu} W. Hurewicz. 
  Ergodic Theorem Without Invariant Measure,
  Annals of Mathematics, Second Series, Vol. 45, No. 1 (Jan., 1944), 192-206.

\bibitem{JT} E. Jarvenpaa, T. Tolonen. 
   Relations between natural and observable measures, 
   Nonlinearity 18:2 (2005), 897-912.

\bibitem{Ka} K. Kaneko. 
    Globally coupled chaos violates the law of large numbers but not the central limit theorem. 
    Phys. Rev. Letters 65 (1990), 1391-1394. 

\bibitem{Ke} G. Keller. 
    An ergodic theoretic approach to mean field coupled maps,
    Progress in Probability, 46 (2000), 183-208.

\bibitem{Si} I.P. Kornfeld,, S.V. Fomin, Ya.G. Sinai.
      Ergodic Theory, Springer (1982). % ISBN 10: 0387905804 ISBN 13: 9780387905808

\bibitem{Mi} M. Misiurewicz. 
  Ergodic natural measures,  
  %"Algebraic and Topological Dynamics" 
  Contemporary Mathematics, 385  
  % eds. S. Kolyada, Y. Manin and T. Ward, Amer. Math. Soc., Providence, RI 
  (2005), 1-6.

\bibitem{OU} J.C. Oxtoby, S.M. Ulam.
       Measure-preserving homeomorphisms and metrical transitivity.
       Ann. of Math. (2), 42 (1941), 874-920.

\bibitem{Ro} V. A. Rohlin, 
  On the fundamental ideas of measure theory, 
  Amer. Math. Soc. Translation 71 (1952), 55 pp;
  (Russian) Mat. Sbornik N.S. 25(67) (1949), 107–150.
%Рохлин, В. А. Об основных понятиях теории меры, Мат. сб., 25:1 (1949), 107–150.

\bibitem{Ru} D. Ruelle. 
     Historical behaviour in smooth dynamical systems, 
     in Global analysis of dynamical systems, edited by 
    H. Broer, B. Krauskopf and G. Vegter, Inst. of Physics, London, 2001, pp.63-66.

%??
%\bibitem{Sch} D. Schnellmann. 
%  Typical points for one-parameter families of piecewise expanding maps of the interval, 
%  Discrete and Continuous Dynamical Systems - Ser A (DCDS-A), 31:3 (2011), 877-911.  

%[FP] B. Faller and C.-E. Pfister. 
%  A point is normal for almost all maps $\beta x + \alpha$ mod 1 or generalized
%  $\beta$-maps, Erg. Theory and Dyn. Systems 29 (2009), 1529-1547.

%??
%\bibitem{WK} C. Wolf, W. Ingle and J.Kaufmann.
%  Natural invariant measures, divergence points and dimension in one-dimensional holomorphic dynamics,  
%  Ergodic Theory and Dynam. Systems 29:4 (2009), 1235-1255.

\end{thebibliography}
\end{document}